\documentclass[reqno]{amsart}
\usepackage{amssymb}

\newtheorem{theorem}{Theorem}

\newtheorem{corollary}[theorem]{Corollary}

\theoremstyle{definition}

\theoremstyle{remark}
\newtheorem{remark}{Remark}

\begin{document}

\title[A sum of theta functions]
{A multiparameter summation formula
\\ for  Riemann  theta functions}

\author{Vyacheslav P. Spiridonov}
 \address{Bogoliubov Laboratory of Theoretical Physics, Joint
Institute for Nuclear Research, Dubna, Moscow Region 141980, Russia}

 \thanks{{\em Date:} November 2004; Proceedings of the Workshop on Jack, Hall-Littlewood
  and Macdonald polynomials (Edinburgh, September 23-26, 2003), Contemp. Math.,
 to appear.}

\copyrightinfo{2004}{V.P. Spiridonov}

\begin{abstract}
We generalize Warnaar's elliptic extension of a Macdonald multiparameter
summation formula to Riemann surfaces of arbitrary genus.
\end{abstract}

\maketitle

We start by a brief outline of the general theory of
hypergeometric type series built from Jacobi theta functions \cite{S1,S2}
(its extension to integrals \cite{S4}
is not touched at all). Within this approach,
univariate elliptic hypergeometric series are defined as formal series
$\sum_{n}c_n$ for which $h(n)=c_{n+1}/c_n$ is an elliptic function
of $n$ considered as a continuous complex variable.
Normalizing $c_0=1$, we see that all coefficients $c_n$ are obtained
as products of $h(k)$ or $1/h(k)$ for different $k\in\mathbb{Z}$.

Any elliptic function of order $r+1$ can be represented as \cite{WW}:
\begin{eqnarray}
&& h(n)=z\frac{\theta_1(u_0+n,\ldots,u_r+n;\sigma,\tau)}
{\theta_1(v_0+n,\ldots,v_r+n;\sigma, \tau)}=z
\frac{\theta(t_0q^n,\ldots,t_rq^n ;p)}
{\theta(w_0q^n,\ldots,w_rq^n;p)},
\label{n=1} \\
&& \theta_1(u_0,\ldots, u_k;\sigma,\tau)=
\prod_{i=0}^k\theta_1(u_i;\sigma,\tau),
\qquad
\theta(t_0,\ldots, t_k;p)=\prod_{i=0}^k\theta(t_i;p),
\nonumber\end{eqnarray}
with a free variable $z\in\mathbb{C}$, the base $p=e^{2\pi i\tau }$
satisfying the constraint $|p|<1$ (or $\text{Im}(\tau)>0$), and
arbitrary $q=e^{2\pi i\sigma}$. The function
\begin{eqnarray} \nonumber
\theta_1(u;\sigma,\tau) &=&
-i\sum_{k=-\infty}^\infty (-1)^kp^{(2k+1)^2/8}q^{(k+1/2)u}
\\ \nonumber
&=& ip^{1/8} q^{-u/2}\: (p;p)_\infty\: \theta(q^{u};p), \quad
u\in\mathbb{C},
\label{theta1}\end{eqnarray}
is the standard Jacobi $\theta_1$-function.
The short theta function $\theta(a;p)$ has the form
\begin{equation}\label{j-theta}
\theta(a;p)=(a;p)_\infty(pa^{-1};p)_\infty,
\qquad (a;p)_\infty=\prod_{n=0}^\infty(1-ap^n),
\end{equation}
and obeys the properties
\begin{equation}\label{fun-rel}
\theta(pa;p)=\theta(a^{-1};p)=-a^{-1}\theta(a;p).
\end{equation}
We use also the convenient elliptic numbers notation:
$$
[u;\sigma,\tau]\equiv \theta_1(u;\sigma,\tau), \quad
\text{or}\quad [u]\equiv \theta_1(u).
$$

The parameters $u_i$ and $v_i$ in (\ref{n=1}) are defined modulo
shifts by $\sigma^{-1}$, which will not be indicated in the formulas
below. They are connected to $t_i$ and $w_i$ as
$$
t_i=q^{u_i}, \qquad w_i=q^{v_i}
$$
and satisfy the constraint
\begin{equation}
\sum_{i=0}^r(u_i-v_i)=0, \quad \text{or}\quad
\prod_{i=0}^rt_i=\prod_{i=0}^rw_i,
\label{balance}\end{equation}
which guarantees double periodicity of the meromorphic function $h(n)$:
$$
h(n+\sigma^{-1})=h(n),\qquad
h(n+\tau\sigma^{-1})=h(n).
$$

For the entire function $[u]$ we have $[-u]=-[u]$,
\begin{equation}
[u+\sigma^{-1}] = -[u], \qquad
[u+\tau\sigma^{-1}] = -e^{-\pi i\tau-2\pi i\sigma u}[u],
\label{quasi}\end{equation}
and
\begin{eqnarray}\label{modular-tr1}
&& [u;\sigma,\tau+1] = e^{\pi i/4} [u;\sigma,\tau],\qquad
\\  \label{modular-tr2}
&& [u;\sigma/{\tau},-1/\tau]
= -i(-i\tau)^{1/2}e^{\pi i\sigma^2u^2/\tau} [u;\sigma,\tau],
\end{eqnarray}
where the sign of $(-i\tau)^{1/2}$ is fixed from the
positivity of its real part condition.
The latter two transformations generate the modular $PSL(2,\mathbb{Z})$-group,
\begin{equation}\label{sl2}
\tau\to \frac{a\tau+b}{c\tau+d},\qquad \sigma\to \frac{\sigma}{c\tau+d},
\end{equation}
where $a,b,c,d \in\mathbb{Z}$ and $ad-bc=1$.

For the unilateral series $\sum_{n\in\mathbb{N}}c_n$, we conventionally
normalize $w_0=q$ (or $v_0=1$). Then, the elliptic function $h(n)$
generates the single variable {\em elliptic hypergeometric series}:
\begin{equation}
_{r+1}E_r\left({t_0,\ldots, t_{r}\atop w_1,\ldots,w_r};q,p; z\right)
= \sum_{n=0}^\infty \frac{\theta(t_0,t_1,\ldots,t_{r};p;q)_n}
{\theta(q,w_1,\ldots,w_r;p;q)_n}\, z^n.
\label{theta-series}\end{equation}
The elliptic shifted factorials are defined as
$$
\theta(t_0,\ldots,t_k;p;q)_n\equiv \prod_{m=0}^k\prod_{j=0}^{n-1}
\theta(t_mq^j;p),
$$
or in the additive form:
$$
[u_0,\ldots,u_k]_n\equiv \prod_{m=0}^k\prod_{j=0}^{n-1}
[u_m+j;\sigma,\tau].
$$

A further generalization of these functions is constructed from $h(n)$ equal
to an arbitrary meromorphic Jacobi form in the sense of Eichler and Zagier
\cite{EZ}. The corresponding series were called in \cite{S1} as {\em theta
hypergeometric series}. In this classification, theta hypergeometric series
(\ref{theta-series}) are called balanced if their parameters
satisfy the constraint (\ref{balance}).

Elliptic hypergeometric series (\ref{theta-series}) are called
{\em totally elliptic} if  $h(n)$ is also an elliptic function of all free
parameters entering its $\theta_1$-functions arguments. It is not possible
to have this property if all $u_i$ and $v_i$ are independent. Therefore
such a requirement results in the
following additional constraints \cite{S1}:
\begin{equation}\label{well-poised}
u_0+1=u_1+v_1=\ldots=u_{r}+v_r, \quad \text{or}\quad
qt_0=t_1w_1=\ldots=t_rw_r.
\end{equation}
These relations are known as the {\em well-poisedness} condition
for plain and basic hypergeometric series \cite{GR}.
Thus the total ellipticity concept sheds some light
on the origin of the restrictions (\ref{well-poised}).
Totally elliptic hypergeometric series are automatically modular invariant.

The balancing condition for well-poised series contains a sign ambiguity:
$t_1\cdots t_r$ $= \pm q^{(r+1)/2}t_0^{(r-1)/2}$. It clearly shows also
distinguished character of the parameter $t_0$. It is convenient to express $t_r$
in terms of other parameters. In this case the function $h(n)$ is invariant
with respect to the independent $p$-shifts $t_0\to p^2t_0$ and $t_j\to pt_j,\,
j=1,\ldots,r-1$. However, for odd $r=2m+1$ we can reach the $t_0\to pt_0$
invariance provided we resolve the ambiguity in the balancing condition in
favor of the form $t_1\cdots t_{2m+1}= + q^{m+1}t_0^{m}$, which is precisely
the condition necessary for obtaining nontrivial hypergeometric identities.
Thus the notion of total ellipticity (with equal periods) uniquely
distinguishes the correct form of the balancing condition.

The elliptic analog of the {\em very-well-poisedness} condition
consists in adding to (\ref{well-poised}) of four constraints \cite{S1}:
\begin{equation}
 t_{r-3}=t_0^{1/2}q,\quad t_{r-2}=-t_0^{1/2}q, \quad
 t_{r-1}=t_0^{1/2}qp^{-1/2},\quad t_{r}=-t_0^{1/2}qp^{1/2}.
\label{vwp}\end{equation}
After simplifications, the very-well-poised theta hypergeometric series
$_{r+1}E_r$ take the form (no balancing condition is assumed)
\begin{equation}
_{r+1}E_r(\ldots)
= \sum_{n=0}^\infty \frac{\theta(t_0q^{2n};p)}{\theta(t_0;p)}
\prod_{m=0}^{r-4}\frac{\theta(t_m;p;q)_n}
{\theta(qt_0/t_m;p;q)_n}\, (-qz)^n.
\label{vwp-1}\end{equation}
In the limit $p\to 0$, functions (\ref{vwp-1}) are reduced to the
very-well-poised  basic hypergeometric series $_{r-1}\varphi_{r-2}$
of the argument $-qz$ \cite{GR}. For even $r$ the balancing condition for
(\ref{vwp-1}) has the form $t_1\cdots t_{r-4}= \pm q^{(r-7)/2}t_0^{(r-5)/2}$.
For odd $r=2m+1$, it is appropriate to call (\ref{vwp-1}) balanced if
$$
\prod_{j=1}^{2m-3}t_j= q^{m-3}t_0^{m-2}.
$$
Then the function $h(n)$ is elliptic in $u_0$, i.e. it
is invariant with respect to the shift $t_0\to pt_0$.
Conditions \eqref{vwp} do not spoil this property
because the change $t_0\to pt_0$ is equivalent (due to the permutational
invariance) to a replacement of parameters \eqref{vwp} by
$t_{r-3},pt_{r-2},pt_{r-1},t_r$ and, so, it does not have an effect
upon $h(n)$. The modular invariance (which is not immediately evident
due to the dependence of parameters in \eqref{vwp} on the base $p$)
is preserved due to a similar reasoning.

For the first time, the very-well-poised elliptic hypergeometric series
(with a different way of counting the number of parameters and the choice
$z=-1$) have been
considered by Frenkel and Turaev \cite{FT}. Their work was inspired
by exactly solvable statistical mechanics models built by Date et al
\cite{D-O}. In an independent setting, such functions have been derived
by Zhedanov and the author by solving a three term recurrence relation
for a self-similar family of biorthogonal rational functions \cite{SZ}.

As shown by Frenkel and Turaev \cite{FT}, the following summation formula is true:
\begin{eqnarray} \nonumber
\lefteqn{\sum_{k=0}^n
\frac{\theta (t_0q^{2k};p)}{\theta(t_0;p)}
\prod_{m=0}^5 \frac{\theta (t_m ;p;q)_k}
     {\theta (qt_0 t_m^{-1};p;q)_k}\, q^k }
&& \\ && \makebox[4em]{}
=\frac{\theta (qt_0;p;q)_n\prod_{1\leq r<s\leq 3}
      \theta (qt_0/t_rt_s;p;q)_n}
{\theta (qt_0/t_1t_2t_3;p;q)_n
\prod_{r=1}^3\theta (qt_0/t_r;p;q)_n},
\label{ft-sum}\end{eqnarray}
where $\prod_{i=1}^5 t_i=qt_0^2$ and $t_4=q^{-n}$, $n\in \mathbb{N}$.
This means that the special terminating very-well-poised balanced $_{10}E_9$
series with $z=-1$ is summable.
If we take the limit $t_3\to t_0/t_1t_2$, then the termination condition
cancels out and \eqref{ft-sum} becomes a summation formula for an
indeterminate $_8E_7$ series
\begin{equation}
\sum_{k=0}^n\frac{\theta (t_0q^{2k};p)}{\theta(t_0;p)}
\prod_{m=0}^3\frac{\theta (t_m;p;q)_k}
     {\theta (qt_0/t_m;p;q)_k}\, q^k
=\frac{\theta (qt_0,qt_1,qt_2,qt_0/t_1t_2;p;q)_n}
{\theta(q,qt_0/t_1,qt_0/t_2,qt_1t_2;p;q)_n}.
\label{8E7}\end{equation}
For $p\to 0$, we obtain an indeterminate summation formula for
a very-well-poised balanced $_6\varphi_5$ basic hypergeometric series.
In the additive notation, equality \eqref{8E7} takes the form
\begin{eqnarray}\nonumber
&& \sum_{k=0}^n\frac{[u_0+2k]}{[u_0]}
\frac{[u_0,u_1,u_2,u_0-u_1-u_2]_k}{[1,1+u_0-u_1,1+u_0-u_2, 1+u_1+u_2]_k}
\\ && \makebox[4em]{}
=\frac{[1+u_0,1+u_1,1+u_2,1+u_0-u_1-u_2]_n}
{[1,1+u_0-u_1,1+u_0-u_2,1+u_1+u_2]_n}.
\label{8E7-add}\end{eqnarray}

The general theory of multiple elliptic
hypergeometric series is much less developed. Several examples of such
series associated with root systems have been investigated by
Warnaar \cite{Wa}, van Diejen and the author \cite{DS1,S1,S3},
 Rosengren \cite{Ro}, Kajihara and Noumi \cite{KN}, Rains \cite{Ra},
with the latest contributions coming from the Gustafson's
talk at the present conference.

We expect that the approach of \cite{S1} to elliptic
hypergeometric series is generalizable to
multidimensional Riemann theta functions associated with algebraic curves.
In this case, we can take as $h(n)$ a meromorphic
function upon a compact Riemann surface $S$ of an arbitrary genus $g$.
However, the things are much less transparent in this situation.
In \cite{S3}, we started to discuss such generalizations and
suggested an extension of sum \eqref{8E7-add} to a series built from the
Riemann theta functions. However, there was a flaw in the
formulation of this result. The main purpose of the present note is
to lift the Macdonald multiparameter sum described in Theorem 2.27 of
\cite{BM} to Riemann surfaces of arbitrary
genus and to correct a related statement of \cite{S3}.

The Riemann theta function of $g$ variables $u_1,\ldots,u_g\in \mathbb{C}$
with characteristics $\mathbf{\alpha}=(\alpha_1,\ldots,\alpha_g)$ and
$\mathbf{\beta}=(\beta_1,\ldots,$ $\beta_g)$ is defined by the $g$-fold
series of the form
\begin{eqnarray}\nonumber 
&& \Theta_{\mathbf{\alpha},\mathbf{\beta}}(\mathbf{u};\Omega)
=\sum_{\mathbf{n}\in \mathbb{Z}^g}
\exp \{ \pi i\sum_{j,k=1}^g (n_j+\alpha_j)\Omega_{jk}(n_k+\alpha_k)
\\ && \makebox[10em]{}
+2\pi i\sum_{j=1}^g(u_j+\beta_j)(n_j+\alpha_j) \},
\label{theta}\end{eqnarray} 
where $\Omega_{jk}$ is a symmetric matrix of periods.
For general abelian varieties the matrix $\Omega_{jk}$ is arbitrary,
and for the Jacobians associated with the Riemann surfaces,
it is generated by a basis of holomorphic
differentials ${\bf \omega}= (\omega_1,\ldots,\omega_g)$
(see, e.g., \cite{Mu}). It is evident that for arbitrary $j$ we have
\begin{equation}
\Theta_{\mathbf{\alpha},\mathbf{\beta}}(u_1,\ldots,u_j+1,\ldots, u_g;\Omega)=
e^{2\pi i \alpha_j}\Theta_{\mathbf{\alpha},\mathbf{\beta}}(\mathbf{u};\Omega).
\label{per-1}\end{equation}
Analogously, we have
\begin{equation}
\Theta_{\mathbf{\alpha},\mathbf{\beta}}(u_1+\Omega_{1k},\ldots, u_g
+\Omega_{gk};\Omega)=
e^{-\pi i\Omega_{kk}-2\pi i(\beta_k+u_k)}\Theta_{\mathbf{\alpha},\mathbf{\beta}}
(\mathbf{u};\Omega),
\label{per-2}\end{equation}
where $k=1,\ldots, g$.

We denote as $\Gamma_{1,2}$ a subgroup of the symplectic modular group
$Sp(2g,\mathbb{Z})$ generated by the matrices
$$
\gamma=\left(\begin{array}{cc} a & b \\ c & d \end{array}\right)
\in Sp(2g,\mathbb{Z}),
$$
such that $diag(a^tb)=diag(c^td)=0\,\text{mod}\, 2$.
The action of this group upon
the matrix of periods $\Omega$ and the arguments $\mathbf{u}$ of the
theta function is defined as
\begin{equation}
\Omega'=(a\Omega+b)(c\Omega+d)^{-1},\qquad
\mathbf{u}'=\mathbf{u}^t(a\Omega+b)^{-1}.
\label{sp}\end{equation}
Analogously, we define the tranformed characteristics
$$
\left(\begin{array}{c} \mathbf{\alpha}' \\ \mathbf{\beta}' \end{array}\right)
= \left(\begin{array}{cc} d & -c \\ -b & a \end{array}\right)
\left(\begin{array}{c} \mathbf{\alpha} \\ \mathbf{\beta} \end{array}\right)
+\frac{1}{2}\left(\begin{array}{c} diag(c^td) \\ diag(a^td) \end{array}\right).
$$
Then, the $\Gamma_{1,2}$ group transformation law for theta functions has
the form
\begin{equation}
\Theta_{\mathbf{\alpha}',\mathbf{\beta}'}(\mathbf{u}';\Omega')
=\zeta\sqrt{\det(c\Omega+d)}e^{\pi i \mathbf{u}^t(c\Omega+d)^{-1}c\mathbf{u}}
\Theta_{\mathbf{\alpha},\mathbf{\beta}}(\mathbf{u};\Omega),
\label{sp-mod}\end{equation}
where $\zeta$ is an eighth root of unity \cite{Mu}.

We denote as
\begin{equation}
v_j(a,b)\equiv \int_a^b\omega_j,\qquad a,b\in S,
\label{ab-int}\end{equation}
abelian integrals of the first kind. The characteristics $\alpha_j,
\beta_j\in (0, 1/2)$, such that $4\sum_{j=1}^g\alpha_j\beta_j=1\,
(\text{mod}\; 2)$, are called odd. We denote theta functions with
arbitrary (nonsingular) odd characteristics as $[{\bf u}],\; {\bf u}\in \mathbb{C}^g,$
and take the convention that
$
[\mathbf{u}_1,\ldots,\mathbf{u}_k]= \prod_{j=1}^k[\mathbf{u}_j].
$
For such functions, we have
$[-{\bf u}]=-[{\bf u}]$. It is known that theta functions on Riemann surfaces
satisfy the Fay identity \cite{F}:
\begin{eqnarray} \nonumber
&& [{\bf u}+\mathbf{v}(a,c),{\bf u}+\mathbf{v}(b,d),
\mathbf{v}(c,b),\mathbf{v}(a,d)]
\\ \label{Fay} && \makebox[2em]{}
+[{\bf u}+\mathbf{v}(b,c),{\bf u}+\mathbf{v}(a,d),
\mathbf{v}(a,c),\mathbf{v}(b,d)]
\\ && \makebox[4em]{}
= [{\bf u},{\bf u}+\mathbf{v}(a,c)+\mathbf{v}(b,d),
\mathbf{v}(c,d),\mathbf{v}(a,b)],
\nonumber\end{eqnarray}
valid for arbitrary ${\bf u}\in \mathbb{C}^g$ and $a,b,c,d\in S$
(note that $\mathbf{v}(a,c)+\mathbf{v}(b,d)=
\mathbf{v}(b,c)+\mathbf{v}(a,d)$). Here we can replace theta functions
independent on the variable $\bf u$ by prime forms $E(x,y)$ since the
cross ratios of these theta functions coincides
with the cross ratios of the corresponding prime forms (and therefore they
are independent on the choice of odd characteristic for theta functions).
However, we find it more convenient to work with one theta function
$[{\bf u}]$. If we consider the function $[{\bf u}+\mathbf{v}(a,b)]$, then,
in the appropriate normalization, transformations \eqref{per-1} and
\eqref{per-2} correspond to cyclic moves of the point $b$ (or $a$)
on the Riemann surface along the $A_i$ and $B_i$ contours, such that
$\int_{A_i}\omega_j=\delta_{ij}$ and
$\int_{B_i}\omega_j=\Omega_{ij}$, respectively \cite{Mu}.

We would like to describe a summation formula for a
particular series built from the Riemann theta functions.
\begin{theorem} We take a nonnegative integer $n$
and consider $n+1$ arbitrary variables $\mathbf{z}_k\in\mathbb{C}^g$ and
$4n+4$ different points on a Riemann surface $a_k, b_k, c_k, d_k\in S$,
$k=0,\ldots, n$. Then the following multiparameter summation formula
for Riemann theta functions on algebraic curves takes place
\begin{eqnarray}\nonumber
&& \sum_{k=0}^n
[\mathbf{z}_k+\mathbf{v}(b_k,c_k),\mathbf{z}_k+\mathbf{v}(a_k,d_k),
\mathbf{v}(a_k,c_k),\mathbf{v}(b_k,d_k)]
\\ \nonumber && \makebox[2em]{} \times
\prod_{j=0}^{k-1}
[\mathbf{z}_j,\mathbf{z}_j+\mathbf{v}(a_j,c_j)+\mathbf{v}(b_j,d_j),
\mathbf{v}(c_j,d_j),\mathbf{v}(a_j,b_j)]
\\ \label{sum} && \makebox[2em]{} \times
\prod_{j=k+1}^{n}
[\mathbf{z}_j+\mathbf{v}(a_j,c_j),\mathbf{z}_j+\mathbf{v}(b_j,d_j),
\mathbf{v}(c_j,b_j),\mathbf{v}(a_j,d_j)]
\\ \nonumber &&
=\prod_{k=0}^n
[\mathbf{z}_k,\mathbf{z}_k+\mathbf{v}(a_k,c_k)+\mathbf{v}(b_k,d_k),
\mathbf{v}(c_k,d_k),\mathbf{v}(a_k,b_k)]
\\ && \makebox[2em]{}
-\prod_{k=0}^n [\mathbf{z}_k+\mathbf{v}(a_k,c_k),
\mathbf{z}_k+\mathbf{v}(b_k,d_k),
\mathbf{v}(c_k,b_k),\mathbf{v}(a_k,d_k)].
 \nonumber\end{eqnarray}
\end{theorem}

\begin{proof}
For proving equality \eqref{sum}, we denote as $f_l^{(n)}$ and $f_r^{(n)}$
its left- and right-hand sides, respectively. We define also
\begin{eqnarray*}
&& g_k=[\mathbf{z}_k,\mathbf{z}_k+\mathbf{v}(a_k,c_k)+\mathbf{v}(b_k,d_k),
\mathbf{v}(c_k,d_k),\mathbf{v}(a_k,b_k)], \\
&& h_k= [\mathbf{z}_k+\mathbf{v}(a_k,c_k),
\mathbf{z}_k+\mathbf{v}(b_k,d_k),
\mathbf{v}(c_k,b_k),\mathbf{v}(a_k,d_k)],
\end{eqnarray*}
so that $f_r^{(n)}=\prod_{k=0}^n g_k -\prod_{k=0}^n h_k$.

For $n=0$, equality \eqref{sum} is reduced to the identity \eqref{Fay}.
Suppose that \eqref{sum} is true for $n=0,\ldots, N-1$,
$N\geq 1$. Then we have by induction
\begin{eqnarray*}
&& f_l^{(N)}=h_Nf_l^{(N-1)}
\\ && \makebox[2em]{}
+[\mathbf{z}_N+\mathbf{v}(b_N,c_N),\mathbf{z}_N+\mathbf{v}(a_N,d_N),
\mathbf{v}(a_N,c_N),\mathbf{v}(b_N,d_N)]
\prod_{j=0}^{N-1}g_j
\\  && \makebox[4em]{}
=\xi_N\prod_{k=0}^{N-1}g_k -\prod_{k=0}^N h_k,
\end{eqnarray*}
where
$$
\xi_N=h_N + [\mathbf{z}_N+\mathbf{v}(b_N,c_N),
\mathbf{z}_N+\mathbf{v}(a_N,d_N),
\mathbf{v}(a_N,c_N),\mathbf{v}(b_N,d_N)].
$$
Using the Fay identity \eqref{Fay}, we find that $\xi_N=g_N$.
Therefore, $f_l^{(N)}=\prod_{k=0}^Ng_k-\prod_{k=0}^Nh_k=f_r^{(N)}$,
that is formula \eqref{sum} is valid for arbitrary $n$.
\end{proof}

\begin{remark}
As noted by the referee, this theorem is a special case of the general construction
for telescoping sums:
\begin{equation}
\sum_{k=0}^{n} (x_k-y_k) \prod_{j=0}^{k-1} x_j
\prod_{j=k+1}^{n} y_j = \prod_{j=0}^{n} x_j - \prod_{j=0}^{n} y_j,
\end{equation}
whose proof follows similar lines:
$$
\sum_{k=0}^{n} (x_k-y_k) \prod_{j=0}^{k-1} x_j \prod_{j=k+1}^{n} y_j
 = \sum_{k=0}^{n} \prod_{j=0}^{k} x_j \prod_{j=k+1}^{n} y_j
 - \sum_{k=0}^{n} \prod_{j=0}^{k-1} x_j \prod_{j=k}^{n} y_j
$$
with an obvious cancellation of extra terms in the latter sums.
\end{remark}

For elliptic curves (i.e., for $g=1$) formula \eqref{sum} was
proven by Warnaar \cite{Wa}. Its further degeneration to the
trigonometric level leads to a Macdonald identity, which was
published for the first time in the Bhatnagar-Milne paper
\cite{BM} and which generalizes relations
obtained by Chu in \cite{C}. As shown in \cite{BM,Wa},
equalities of such type with special choices of parameters can be
used for derivations of more fine structured summation and transformation
formulas for bibasic and elliptic hypergeometric series.

We suppose that for some $N>0$ the points $a_N,b_N,c_N,d_N\in S$
are such that we hit a zero of a theta function: $[\mathbf{v}(a_N,b_N)]=0$
or $[\mathbf{v}(c_N,d_N)]=0$. Then, using the antisymmetry
$[\mathbf{v}(c_k,b_k)]=-[\mathbf{v}(b_k,c_k)]$
equality \eqref{sum} can be rewritten as
\begin{eqnarray}\nonumber
&& \sum_{k=0}^N
\frac{[\mathbf{z}_k+\mathbf{v}(b_k,c_k),\mathbf{z}_k+\mathbf{v}(a_k,d_k),
\mathbf{v}(a_k,c_k),\mathbf{v}(b_k,d_k)]}
{[\mathbf{z}_k+\mathbf{v}(a_k,c_k),
\mathbf{z}_k+\mathbf{v}(b_k,d_k),
\mathbf{v}(b_k,c_k),\mathbf{v}(a_k,d_k)]}
\\ && \makebox[2em]{} \times
\prod_{j=0}^{k-1}
\frac{[\mathbf{z}_j,\mathbf{z}_j+\mathbf{v}(a_j,c_j)+\mathbf{v}(b_j,d_j),
\mathbf{v}(c_j,d_j),\mathbf{v}(a_j,b_j)]}
{[\mathbf{z}_j+\mathbf{v}(a_j,c_j),
\mathbf{z}_j+\mathbf{v}(b_j,d_j),
\mathbf{v}(c_j,b_j),\mathbf{v}(a_j,d_j)]}
= 1.
\label{sum2}\end{eqnarray}

Equivalently, for $[\mathbf{v}(c_0,b_0)]=0$ or $[\mathbf{v}(a_0,d_0)]=0$
we obtain the sum
\begin{eqnarray} \nonumber
\lefteqn{ \sum_{k=0}^n
\frac{
[\mathbf{z}_k+\mathbf{v}(b_k,c_k),\mathbf{z}_k+\mathbf{v}(a_k,d_k),
\mathbf{v}(a_k,c_k),\mathbf{v}(b_k,d_k)]}
{[\mathbf{z}_0,\mathbf{z}_0+\mathbf{v}(a_0,c_0)+\mathbf{v}(b_0,d_0),
\mathbf{v}(c_0,d_0),\mathbf{v}(a_0,b_0)]}
} && \\ \nonumber &&  \times
\prod_{j=0}^{k-1}\frac{
[\mathbf{z}_j,\mathbf{z}_j+\mathbf{v}(a_j,c_j)+\mathbf{v}(b_j,d_j),
\mathbf{v}(c_j,d_j),\mathbf{v}(a_j,b_j)]}
{[\mathbf{z}_{j+1}+\mathbf{v}(a_{j+1},c_{j+1}),
\mathbf{z}_{j+1}+\mathbf{v}(b_{j+1},d_{j+1}),
\mathbf{v}(c_{j+1},b_{j+1}),\mathbf{v}(a_{j+1},d_{j+1})]}
\\ && \makebox[4em]{}
=\prod_{k=1}^n\frac{
[\mathbf{z}_k,\mathbf{z}_k+\mathbf{v}(a_k,c_k)+\mathbf{v}(b_k,d_k),
\mathbf{v}(c_k,d_k),\mathbf{v}(a_k,b_k)]}
{[\mathbf{z}_{k}+\mathbf{v}(a_{k},c_{k}),
\mathbf{z}_{k}+\mathbf{v}(b_{k},d_{k}),
\mathbf{v}(c_{k},b_{k}),\mathbf{v}(a_{k},d_{k})]}.
\label{sum3}\end{eqnarray}

In the elliptic case and its further degenerations, such equalities are
useful in searching matrices whose inversions are given
by simple analytical expressions. It is natural
to expect that our relations will lead to a generalization
of some of the Warnaar's results on multibasic hypergeometric sums
and matrix inversions \cite{Wa}. Moreover, it is worth to analyze
possible extensions of the Krattenthaler's determinant
from the Warnaar's elliptic case \cite{Wa} to arbitrary Riemann surfaces.

By special choices of parameters we can try to give to the derived sums
hypergeometric type forms. We consider only two particular cases.
Substituting
$$
\mathbf{z}_k=\mathbf{u}_0+\mathbf{v}(x_k), \quad
b_k=c_0\equiv x_0,\quad c_k\equiv x_k,\quad d_k=d_0,
$$
where $\mathbf{u}_0\in\mathbb{C}^g$,
$\mathbf{v}(x_k)\equiv\mathbf{v}(x_0,x_k)$, $k=0,1,\ldots,$ and
$\mathbf{u}_2\equiv\mathbf{v}(d_0,x_0),$ $\mathbf{u}_1^{(k)}\equiv\mathbf{v}(x_k,a_k)$
into (\ref{sum3}), we obtain the following identity.

\begin{corollary}
\begin{eqnarray}\nonumber
\lefteqn{
\sum_{k=0}^n\frac{[{\bf u}_0+2\mathbf{v}(x_k)]}{[\mathbf{u}_0]}
\frac{[\mathbf{u}_0-\mathbf{u}_1^{(k)}-\mathbf{u}_2,\mathbf{u}_1^{(k)}]}
{[\mathbf{u}_0-\mathbf{u}_1^{(0)}-\mathbf{u}_2,\mathbf{u}_1^{(0)}]}
\prod_{j=0}^{k-1}\Biggl(\frac{[{\bf u}_0+\mathbf{v}(x_j)]}{[\mathbf{v}(x_{j+1})]} }
&& \\ \nonumber && \times
\frac{[{\bf u}_0-{\bf u}_1^{(j)}-{\bf u}_2+\mathbf{v}(x_j),
{\bf u}_1^{(j)}+\mathbf{v}(x_j),{\bf u}_2+\mathbf{v}(x_j)]}
{[{\bf u}_1^{(j+1)}+{\bf u}_2+\mathbf{v}(x_{j+1}),
{\bf u}_0-{\bf u}_1^{(j+1)}+\mathbf{v}(x_{j+1}),
{\bf u}_0-{\bf u}_2+\mathbf{v}(x_{j+1})]}\Biggr)
\\ &&
= \prod_{k=1}^{n}\frac
{[{\bf u}_0+\mathbf{v}(x_k),
{\bf u}_0-{\bf u}_1^{(k)}-{\bf u}_2+\mathbf{v}(x_k),
{\bf u}_1^{(k)}+\mathbf{v}(x_k),{\bf u}_2+\mathbf{v}(x_k)]}
{[\mathbf{v}(x_k),{\bf u}_1^{(k)}+{\bf u}_2+\mathbf{v}(x_k),
{\bf u}_0-{\bf u}_1^{(k)}+\mathbf{v}(x_k),
{\bf u}_0-{\bf u}_2+\mathbf{v}(x_k)]}.
\label{theta-sum}\end{eqnarray}
\end{corollary}

If there would exist a sequence of points $a_k$ such that
$\mathbf{u}_1^{(k)}=\mathbf{u}_1^{(0)}$, then four theta functions
of the second factor in the left-hand side of ({\ref{theta-sum}) would
cancel each other
and we would obtain an exact $g>1$ analog  of the sum (\ref{8E7-add}).
Such a condition was assumed implicitly in the derivation of
the corresponding formula in \cite{S3}.
However, in general this is possible only for $g=1$ (with the choice
$\omega= du, \, a_k=k+a_0, x_k=k+x_0$). As a result, our $g>1$ analog
of (\ref{8E7-add}) (which corrects the related formula of \cite{S3})
does not obey all structural properties of the $g=1$
summation formula. Still, we draw attention to the fact that the
right-hand side of ({\ref{theta-sum}) satisfies analogs of
the balancing and well-poisedness conditions: 1) the sums of arguments of
theta functions in the numerator and denominator are equal to each other;
2) the reciprocal theta functions in the numerator and denominator
have arguments whose sums are equal to ${\bf u}_0+2\mathbf{v}(x_k)$.

Another curious summation formula (inspired by a referee's suggestion)
is obtained after substituting $a_k=a,$ $c_k=b_0=c,$ $b_{k+1}=d_k\equiv x_k$
for $k=0,1,\ldots,$
and ${\bf z}_{k+1}={\bf z}_0+{\bf v}(c,x_k)$ for $k>0$ into (\ref{sum3}).
\begin{corollary}
\begin{eqnarray}\nonumber
\lefteqn{ \sum_{k=1}^{n}
\frac{[{\bf z_0},{\bf z}_0+{\bf v}(c,x_{k-1})+{\bf v}(a,x_k),
{\bf v}(a,c),{\bf v}(x_{k-1},x_k)]}
{[{\bf z}_0+{\bf v}(c,x_{k-1}),{\bf z}_0+{\bf v}(c,x_k),
{\bf v}(a,x_{k-1}),{\bf v}(a,x_k)]} } &&
\\ &&
 = \frac{[{\bf z}_0+{\bf v}(a,x_n),{\bf v}(c,x_n)]}
{[{\bf z}_0+{\bf v}(c,x_n),{\bf v}(a,x_n)]}
 - \frac{[{\bf z}_0+{\bf v}(a,x_0),{\bf v}(c,x_0)]}
{[{\bf z}_0+{\bf v}(c,x_0),{\bf v}(a,x_0)]}
\\   \nonumber &&
 = \frac{[{\bf z}_0,{\bf z}_0+{\bf v}(c,x_0)+{\bf v}(a,x_n),
{\bf v}(a,c),{\bf v}(x_0,x_n)]} {[{\bf z}_0+{\bf v}(c,x_0),{\bf z}_0+{\bf v}(c,x_n),
{\bf v}(a,x_0),{\bf v}(a,x_n)]}.
\label{sum4}\end{eqnarray}
\end{corollary}

It is easy also to see that all the derived sums represent ``totally abelian"
functions, that is they are invariant under arbitrary moves of points on
the Riemann surface along the cycles and appropriate $2g$ shifts of
the variables $\mathbf{z}_k$ (or $\mathbf{u}_0$). Analogously, the
$Sp(2g,\mathbb{Z})$ modular group invariance is evident by construction
(due to the Fay identity properties). As it is clear from our preliminary
analysis, theta hypergeometric series for Riemann surfaces with the genus $g>1$
should obey some principally new features with respect to the elliptic case
and their determination requires some additional effort.

\smallskip

The author is deeply indebted to the organizers of the Workshop on
Jack, Hall-Littlewood and Macdonald polynomials  (Edinburgh,
September 23-26, 2003) for invitation to speak there and to
the referee for pointing to a flaw in the original formulation of
Corollary 2.
This work is supported in part by the Russian Foundation for Basic
Research (grant No. 03-01-00780). Its final version was prepared
during the author's stay at the Max-Planck-Institut f\"ur Mathematik
in Bonn whose hospitality is gratefully acknowledged.


\begin{thebibliography}{000000}


\bibitem{BM} G. Bhatnagar and S. C. Milne, {\em Generalized bibasic
hypergeometric series and their $U(n)$ extensions}, Adv. Math. {\bf 131}
(1997), 188--252.

\bibitem{C} W. C. Chu, {\em Inversion techniques and combinatorial
identities}, Bull. Un. Mat. Ital. {\bf 7} (1993), 737--760.

\bibitem{D-O} E. Date, M. Jimbo, A. Kuniba, T. Miwa,
and M. Okado, {\em Exactly solvable SOS models: local height
probabilities and theta function identities}, Nucl. Phys. B
{\bf 290} (1987), 231--273.

\bibitem{DS1} J. F. van Diejen and V. P. Spiridonov,
 {\em Modular hypergeometric
residue sums of elliptic Selberg integrals},
Lett. Math. Phys. {\bf 58} (2001), 223--238.

\bibitem{EZ} M. Eichler and D. Zagier,
{\em The Theory of Jacobi Forms}, Progress in Math. {\bf 55},
Birkh\"auser, Boston, 1985.

\bibitem{F} J. F. Fay, {\em Theta functions on Riemann surfaces},
Lect. Notes in Math. {\bf 353}, Springer-Verlag, Berlin, 1973.

\bibitem{FT} I. B. Frenkel and V. G. Turaev,
{\em Elliptic solutions of the Yang-Baxter equation and modular
hypergeometric functions}, The Arnold-Gelfand Mathematical
Seminars, Birkh\"auser, Boston, 1997, pp. 171--204.

\bibitem{GR} G. Gasper and M. Rahman,
{\em Basic Hypergeometric Series},
Encyclopedia of Mathematics and its Applications {\bf 35}, Cambridge
Univ. Press, Cambridge, 1990.

\bibitem{KN} Y. Kajihara and M. Noumi, {\em Multiple elliptic
hypergeometric series. An approach from the Cauchy determinant},
Indag. Math. {\bf 14} (2003), 395--421.

\bibitem{Mu} D. Mumford, {\em Tata Lectures on Theta I, II},
Progress in Math. {\bf 28, 43}, Birkh\"auser, Boston, 1983, 1984.

\bibitem{Ra} E. M. Rains, {\em Transformations of
elliptic hypergeometric integrals}, preprint (2003), arXiv:math.QA/0309252.

\bibitem{Ro} H. Rosengren, {\em Elliptic hypergeometric series on root
systems}, Adv. Math. {\bf 181} (2004), 417--447.

\bibitem{S1} V. P. Spiridonov, {\em Theta hypergeometric series},
Asymptotic Combinatorics with Application to Mathematical Physics
(St. Petersburg, July 9--22, 2001), Kluwer, 2002, pp. 307--327.

\bibitem{S2} \bysame, {\em An elliptic incarnation of the Bailey chain},
Internat. Math. Res. Notices, no. 37 (2002), 1945--1977.

\bibitem{S3} \bysame, {\em Modularity and total ellipticity of some
 multiple series of hypergeometric type},
 Theor. Math. Phys. {\bf 135} (2003), 836--848.

\bibitem{S4} \bysame, {\em Theta hypergeometric integrals}, Algebra i
Analiz {\bf 15} (6) (2003), 161--215 (St. Petersburg Math. J.
{\bf 15} (6) (2004), 929--967).

\bibitem{SZ} V. P. Spiridonov and A. S. Zhedanov,
{\em Spectral transformation chains and some new biorthogonal
rational functions}, Commun. Math. Phys. {\bf 210} (2000), 49--83.

\bibitem{Wa} S. O. Warnaar, {\em Summation and transformation formulas for
elliptic hypergeometric series}, Constr. Approx. {\bf 18} (2002), 479--502.

\bibitem{WW} E. T. Whittaker and G. N. Watson, {\em A Course
of Modern Analysis}, Cambridge Univ. Press, Cambridge, 1986.

\end{thebibliography}
\end{document}